\documentclass[final%
,nohyperref%
]{dmtcs-episciences}

\usepackage[utf8]{inputenc}
\usepackage{subfigure}
\usepackage{amssymb}
\usepackage{amsmath}
\usepackage{amsfonts,epsf,amsmath,tikz}
\usetikzlibrary{decorations.pathreplacing,automata,calc,positioning}
\usepackage{framed}

\newcommand{\Mod}[1]{\ (\text{mod}\ #1)}

\newtheorem{theorem}{\bf Theorem}[section]
\newtheorem{corollary}[theorem]{\bf Corollary}
\newtheorem{lemma}[theorem]{\bf Lemma}
\newtheorem{proposition}[theorem]{\bf Proposition}
\newtheorem{conjecture}[theorem]{\bf Conjecture}

\newcommand{\Rc}{\mathcal{R}}
\newcommand{\Gc}{\mathcal{G}}
\newcommand{\Oc}{\mathcal{O}}
\newcommand{\Zb}{\mathbb{Z}}
\newcommand{\Nb}{\mathbb{N}}
\newcommand{\cn}{\chi} 
\newcommand{\GDG}{D_\mathcal{G}}
\newcommand{\GDR}{D_\mathcal{R}}

\author{Sylvain Gravier \affiliationmark{1}
  \and Kahina Meslem \affiliationmark{3}
  \and Simon Schmidt \affiliationmark{2}
  \and Souad Slimani \affiliationmark{2,3}}
 
\title{A New Game Invariant of Graphs: the Game Distinguishing Number}
\affiliation{CNRS, Institut Fourier, UMR 5582, S.F.R. Maths \`a Modeler, Universit\'e Grenoble Alpes, France\\
Institut Fourier, S.F.R. Maths \`a Modeler, Universit\'e Grenoble Alpes, France\\
LaROMaD, S.F.R. Maths \`a Modeler, Facult\'e des Math\'ematiques, U.S.T.H.B, Alger, Algeria}

\keywords{distinguishing number, graph automorphism, combinatorial game, hypercube}
\received{2015-10-20}

\accepted{2017-1-30}

\revised{2015-10-22, 2016-8-1, 2016-11-28}
\begin{document}
\publicationdetails{19}{2017}{1}{2}{1303}
\maketitle
\begin{abstract}
 The distinguishing number of a graph $G$ is a symmetry related graph invariant whose study started two decades ago.
The distinguishing number $D(G)$ is the least integer $d$ such that $G$ has a distinguishing $d$-coloring. A distinguishing $d$-coloring is a
coloring $c:V(G)\rightarrow\{1,\cdots,d\}$ invariant only under the trivial automorphism. In this paper, we introduce a game variant of the distinguishing number. The distinguishing game is a game with two players, the Gentle and the Rascal, with antagonist goals. This game is played on a graph $G$ with a set of $d\in\Nb^*$
colors. Alternately, the two players choose a vertex of $G$ and color it with one of the $d$ colors. The game ends when all the vertices have been colored. 
Then the Gentle wins if the coloring is distinguishing and the Rascal wins otherwise. This game leads to define two new invariants for a graph $G$, which are the minimum numbers of colors needed to ensure that the Gentle has a winning strategy, 
depending on who starts. These invariants could be infinite, thus we start by giving sufficient conditions to have infinite game distinguishing numbers. We also show that for graphs with cyclic automorphism group of prime odd order, both game invariants are finite. 
After that, we define a class of graphs, the involutive graphs, for which the game distinguishing number can be quadratically bounded above by the classical distinguishing number. The definition of this class is closely related to imprimitive actions whose blocks have size $2$. 
Then, we apply results on involutive graphs to compute the exact value of these invariants for hypercubes and even cycles.  Finally, we study odd cycles, for which we are able to compute the exact value when their order is not prime. In the prime order case, we give an upper bound of $3$.
\end{abstract}

\section{Background and definition of the game}

The distinguishing number of a graph $G$ is a symmetry related graph invariant. Its study started two decades ago in a work by Albertson and Collins \cite{albertson}.
Given a graph $G$ the distinguishing number $D(G)$ is the least integer $d$ such that $G$ has a distinguishing $d$-coloring. A distinguishing $d$-coloring is a
coloring $c:V(G)\rightarrow\{1,\cdots,d\}$ invariant only under the trivial automorphism. More generally, we say that an automorphism $\sigma$ of a graph $G$ preserves the coloring $c$ or is a color
preserving automorphism, if $c(u)=c(\sigma(u))$ for all $u\in V(G)$. We denote by ${\rm Aut}(G)$, the automorphism group of $G$.
Clearly, for each coloring $c$ of the vertex set of $G$, the set ${\rm Aut}_c(G)=\{\sigma \in {\rm Aut}(G): c\circ\sigma= c\}$ is a subgroup of ${\rm Aut}(G)$. A coloring $c$ is distinguishing if
${\rm Aut}_c(G)$ is trivial.
The group ${\rm Aut}(G)$ acts naturally on the vertex set of $G$ and this action induces a partition of the vertex set $V(G)$ into orbits.
If $H$ is a subgroup of ${\rm Aut}_c(G)$ then the action of $H$ on $V(G)$ is such that each orbit induced by this action is monochromatic. 
In particular, an automorphism $\sigma$ preserves the coloring if and only if the orbits under the action of the subgroup $<\sigma>$ generated by $\sigma$, are all monochromatic.

In the last couple of years the study of this invariant was particularly flourishing. See \cite{klav_complete09,imrich_cartes_power} for the work of Imrich, Jerebic and Klav\v{z}ar on the distinguishing number of the
Cartesian products or \cite{trenk} for an analog of Brook's theorem.  In~\cite{trenk} Collins and Trenk also introduce a distinguishing coloring which must be a proper coloring. 
The distinguishing number is also studied in a more general context than graphs \cite{chan,klwozh-06}. 
Our goal in this paper is to introduce a game variant of this invariant in the spirit of the game chromatic number $\chi_G$ introduced by Brahms in 1981 (see \cite{faigle}) or of the most recent domination game (see \cite{brkl-2010}).  Even if inventing new game invariants could raise a lot of promising and interesting questions, it might seem artificial at first glance. 
To defend our approach, we recall that game invariants have already proved useful to give a new insight on the classical invariant they are related. We cite in particular \cite{klavbujt}, 
where a greedy like strategy used in the study of the domination game is used to improve several upper bounds on the domination number. See also \cite{appl-col} for an application of the coloring game to the graph packing problem. 

The distinguishing game is a game with two players, the Gentle and the Rascal, with antagonist goals. This game is played on a graph $G$ with a set of $d\in\Nb^*$
colors. Alternately, the two players choose a vertex of $G$ and color it with one of the $d$ colors. The game ends when all the vertices have been colored.
Then the Gentle wins if the coloring is distinguishing and the Rascal wins otherwise. 

This game leads to two new invariants for a graph $G$. The {\em $\Gc$-game distinguishing number} $\GDG(G)$ is the minimum of colors needed to ensure that the Gentle has a winning strategy for the game on $G$, assuming he is playing first.
If the Rascal is sure to win whatever the number of colors we allow, then $\GDG(G)=\infty$. Similarly, the {\em $\Rc$-game distinguishing number} $\GDR(G)$ is the minimum number of colors needed to ensure that the Gentle has a winning strategy, assuming the Rascal is playing first. 

In Section \ref{sec 2}, some basic results about the distinguishing games are given. In particular, we study when the distinguishing numbers are finite or not. This question is far from being easy in general. The third section is devoted to introduce the class of involutive graphs, for which we prove that the $\Rc$-game distinguishing number is finite and more precisely, quadratically bounded above by the classical distinguishing number of the graph. 
This class, closely related to imprimitive action of group with complete block system whose blocks have size $2$, contains a variety of graphs such as hypercubes, even cycles, spherical graphs and even graphs \cite{even,spherical}. 
In Section~\ref{sec:hypercube}, the general result on involutive graphs is used to compute the exact values of the game invariants for hypercubes.
\begin{theorem}\label{theo:hypercube}
 Let $Q_n$ be the hypercube of dimension $n\geq2$.
\begin{enumerate}
 \item We have $\GDG(Q_n)=\infty$.
 \item If $n\geq 5$, then $\GDR(Q_n)=2$. Moreover $\GDR(Q_2)=\GDR(Q_3)=3$.
 \item We have $2\leq \GDR(Q_4)\leq 3$.
\end{enumerate}
\end{theorem}   
Finally, in Section~\ref{sec:cycle}, we solve the problem for cycles, except when the number of vertices is prime. For even cycles, it is a straightforward application of a result on involutive graphs, but for odd cycles new ideas are needed.
\begin{theorem}\label{theo:cycle}
 Let $C_n$ be a cycle of order $n\geq 3$.
\begin{enumerate}
\item If $n$ is even, then $\GDG(C_n)=\infty$.  If $n$ is odd, then $\GDR(C_n)=\infty$.
 \item If $n$ is even and $n\geq 8$, then $\GDR(C_n)=2$. Moreover $\GDR(C_4)=\GDR(C_6)=3$.
 \item If $n$ is odd, not prime and $n\geq 9$, then $\GDG(C_n)=2$.
 \item If $n$ is prime and $n\geq 5$, then $\GDG(C_n)\leq 3$. Moreover $\GDG(C_3)=\infty$ and $\GDG(C_5)=\GDG(C_7)=3$.
\end{enumerate}
\end{theorem}
%
%

\section{Basic results}\label{sec 2}
This section is devoted to basic results about the distinguishing games. Especially, we are interested in determining when the $\Gc$-game and $\Rc$-game distinguishing numbers are infinite.

It arises directly from the definition that $D(G)\leq \GDG(G)$ and $D(G)\leq \GDR(G)$, for any graph $G$. Another straightforward remark is that the game distinguishing number of a graph and its complement are equal. 
One first natural question is, if having a winning strategy with $k$ colors is a growing property. In other words, has the Gentle a winning strategy if
$k\geq \GDG$ or $k\geq \GDR$ colors are allowed? The answer is not as obvious as it looks at first glance. We recall that for the game chromatic number $\cn_g$,
it is not known if there is always a winning strategy if we play with $k\geq \cn_g$ colors. However, for the game distinguishing numbers we easily show
they  both have this growing property.
\begin{proposition}
 Let $G$ be a graph and $k$ a positive integer. If $k\geq\GDG(G)$ (resp. $k\geq \GDR (G)$), then the Gentle has a winning strategy with $k$ colors, if he starts (resp. the Rascal starts).
\end{proposition}
\proof In order to win with $k$ colors, the Gentle plays the same winning strategy he would have played with $\GDG(G)$ (resp. $\GDR(G)$) colors, except when the Rascal chooses
a color strictly greater than $\GDG(G)$ (resp. $\GDR(G)$). In that case, the Gentle plays the winning move he would have played if the Rascal had played the color $1$.
Let $c$ be the final coloring in the game with $k$ colors with respect to the above strategy. Let $\tilde{c}$ be the coloring defined for all $v\in V(G)$ as follows:
$$ \tilde c(v)=\begin{cases}
                                                     c(v)  & \text{ if }  c(v)\leq \GDG(G) \text{ (resp. }c(v) \leq \GDR(G))\\
                                                     1 &  \text{otherwise} \\
                                                    \end{cases}.$$

The coloring $\tilde c$ is a distinguishing coloring since it is obtained by following a winning strategy for the game with only $\GDG(G)$
(resp. $\GDR(G)$) colors. Let $\sigma$ be an automorphism which preserves $c$. It is clear that $\sigma$ preserves $\tilde c$ too.
Since $\tilde c$ is distinguishing, $\sigma$ must be the identity.
This shows that $c$ is a distinguishing $k$-coloring and that the Gentle has a winning strategy with $k$ colors.
\qed

As mentioned in the definition, $\GDG$ and $\GDR$ could possibly be infinite. For example, if the graph has an automorphism of order two, the Rascal can win using a strategy 
close to the Tweedledee Tweedledum one for the sum of opposite combinatorial games. He uses the involutive automorphism to copy the move the Gentle just made. 
\begin{proposition} \label{prop:order2}
If $G$ is a graph with a nontrivial  automorphism of order $2$, then:
\begin{enumerate}
   \item $\GDG(G)=\infty$ if $|V(G)|$ is even,
   \item $\GDR(G)=\infty$, if $|V(G)|$ is odd.
\end{enumerate}
\end{proposition}
\proof Let $A$ be the set of vertices fixed by $\sigma$: $A=\{v \in G: \sigma(v)=v\}$. Note that $|V(G)|$ and $|A|$ have the same parity.
We denote by $r_i$ and $s_i$ the $i$-th  vertex played by the Rascal and the Gentle respectively.

First, assume $|A|$ is even and the Gentle starts. The winning strategy for the Rascal is as follows.
If $s_i$ is in $A$, then the Rascal plays another vertex $r_i$ in $A$ and does not pay attention to the color. This is always possible since $|A|$ is even.
Else the Rascal plays such that $r_i=\sigma(s_i)$ and $c(r_i)=c(s_i)$. Since $\sigma$ has order $2$, the vertices outside $A$ can be grouped in pairs $(u,v)$ with $\sigma(u)=v$
and $\sigma(v)=u$. Moreover the Gentle will be the first to play outside $A$. Hence, such a move is always available for the Rascal.

Now suppose $|A|$ is odd and the Rascal begins. His first move is to color a vertex in $A$. The number of uncolored vertices in $A$ are now even and it is the Gentle's
turn. Then the Rascal wins with exactly the same strategy as above.\qed

For example, for cycles $C_n$ or paths $P_n$ on $n$ vertices, we have $\GDG(C_n)=\GDG(P_n)=\infty$, if $n$ is even, and $\GDR(C_n)=\GDR(P_n)=\infty$, if $n$ is odd.
Note that Proposition \ref{prop:order2} tells nothing about the possible values of $\GDG$ (resp. $\GDR$) when the number of fixed points is odd (resp. even).
For cycles the question is rather complicated. We partially solve it in Section \ref{sec:cycle}. But for paths $P_n$, with $n\geq 2$, it is easy to show that $\GDG(P_n)=2$ (resp. $\GDR(P_n)=2$) if $n$ is odd (resp. even). 

There are also graphs with both $\GDG$ and $\GDR$ infinite, for example graphs with two transpositions moving a common vertex, like the complete graphs on more than $3$ vertices.
Finding graphs with both invariants finite, is a bit less obvious. Before stating the result, we recall that for a given finite group $\varGamma$, there are infinitely many non-isomorphic 
graphs whose automorphism group is $\varGamma$.

\begin{proposition}
If $G$ is a graph such that ${\rm Aut}(G)=\Zb/p\Zb$, where $p$ is prime and $p\geq 3$, then $\GDG(G)=\GDR(G)=2$.
\end{proposition}
\proof
Since $|{\rm Aut}(G)|$ is prime all the orbits under the action of the whole group have either size $1$ or size $p$. Let $\Oc_1, \cdots,\Oc_l$, with $l>0$ be the orbits
of size $p$. Since there are at least three distinct vertices in each orbit $\Oc_i$ and two colors are allowed, the Gentle can always play,
 whoever starts, in a manner which ensures that at least one of the $\Oc_i$ is not monochromatic.

Let $\sigma$ be a nontrivial automorphism of $G$. Since $|{\rm Aut}(G)|$ is prime, $\sigma$ generates ${\rm Aut}(G)$ and the orbits under ${\rm Aut}(G)$ and $<\sigma>$ are
the same. The automorphism $\sigma$ cannot preserve the coloring because in that case all the orbits $\Oc_i$ have to be monochromatic.
This shows that the above strategy is winning for the Gentle and that $\GDG(G)=\GDR(G)=2$.
\qed

\section{Involutive graphs}\label{sec:inv-graph}

In this section, we introduce the new class of involutive graphs. This class contains well known graphs such as hypercubes, even cycles, even toroidal grids and even graphs \cite{even,spherical}. For graphs in this class, we are going to prove that their $\Rc$-game distinguishing number is quadratically bounded above by the distinguishing number.

An {\em involutive graph} $G$ is a graph together with an involution, ${\rm Bar}:V(G)\rightarrow V(G)$, which commutes with all automorphisms
and has no fixed point. In others words:\begin{itemize}
                 \item $ \overline{\overline u}=u$ and $\overline u\neq u$ for every $u\in V(G)$\\
                 \item $\sigma(\overline u)=\overline{\sigma(u)}$ for every $\sigma\in {\rm Aut}(G)$ and $u\in V(G)$ .\\
                 \end{itemize}
The set $\{u,\overline u\}$ will be called a {\em block} and $\bar u$ will be called the {\em opposite} of $u$. An important remark is that, since ${\rm Bar}$ commutes with all automorphisms, the image of a block under an automorphism is also a block. The block terminology comes from imprimitive group action theory. Indeed, if we add the condition that the automorphism group of the graph acts transitively, then an involutive graph is a graph such that this action has a complete block system, whose blocks have size $2$.

 
Let $G$ be an involutive graph on $n$ vertices. We define $G'$ as the graph obtained from $G$ by deleting all the edges and putting a new edge between all vertices and their opposites. 
The resulting graph is a disjoint union of $n/2$ copies of $K_2$. Moreover, it is also an involutive graph and it has the same blocks as $G$. The next proposition states that these graphs are in a sense the worst involutive graphs for the Gentle.
\begin{proposition}\label{prop:worst-inv}
 If $G$ is an involutive graph then $D(G)\leq D(G')$ and $\GDR(G)\leq\GDR(G')$.
\end{proposition}
\proof It is enough to prove that any distinguishing coloring $c$ of $G'$ is also a distinguishing coloring of $G$. Let $\sigma$ be a nontrivial automorphism of $G$. Assume first, there are two disjoint blocks $\{u,\bar u\}$ and $\{v,\bar v\}$, with $u,v\in V(G)$, such that $\sigma(\{u,\bar u\})=\{v,\bar v\}$.
Since the coloring $c$ is distinguishing for $G'$, we have $c(\{u,\bar u\})\neq c(\{v,\bar v\})$. Hence, $\sigma$ does not preserve $c$. Suppose now that there are no two such blocks. 
Since $\sigma$ is not trivial, there exists $u\in V(G)$ such that $\sigma(u)=\bar u$. But $u$ and $\bar u$ must have distinct colors, otherwise the transposition which switch them would be a nontrivial color preserving automorphism of $G'$. 
Therefore, $\sigma$ is not a color preserving automorphism. We conclude that $c$ is also a distinguishing coloring for $G$.\qed 

For graphs which are disjoint union of $K_2$, we are able to compute the exact value of $\GDG$ and $\GDR$. 
\begin{proposition}\label{prop:K2}
 If $G$ is the disjoint union of $n\geq1$ copies of $K_2$, then $\GDG(G)=\infty$ and $\GDR(G)=n+1$.
\end{proposition}
\proof Such a graph $G$ has an automorphism of order $2$, with no fixed point, the one which exchanges the two vertices in each copy of $K_2$. By Proposition~\ref{prop:order2}, we directly conclude that $\GDG(G)=\infty$.

We are now going to prove that $\GDR(G)=n+1$. First, assume that only $k<n+1$ colors are allowed during the game. The Rascal's strategy is to play always the color $1$
 in a way that at least one vertex of each copy of $K_2$ is colored with $1$. Since there is strictly less than $n+1$ different colors, whatever the Gentle plays
there will be two distinct $K_2$ colored with the same pair of colors or there will be a monochromatic $K_2$. If there is a monochromatic $K_2$,
 the transposition which permutes its two vertices is clearly a color preserving automorphism. In the other case, there are four distinct vertices $u_1,u_2,v_1,v_2\in V(G)$,
such that $u_1$ and $u_2$ (resp. $v_1$ and $v_2$) are in the same $K_2$. Moreover, $u_1$ and $v_1$ (resp. $u_2$ and $v_2$) share the same color. The automorphism $\sigma$, 
defined by $\sigma(u_i)=v_i$, $\sigma(v_i)=u_i$, with $i\in\{1,2\}$, and $\sigma(w)=w$ for $w\in V(G)\setminus\{u_1,u_2,v_1,v_2\}$ preserves the coloring. Hence, the strategy described is winning for the Rascal. 
In conclusion, $\GDR(G)\geq n+1$.

We now allow $n+1$ colors. We recall that the Rascal plays first. The winning strategy for the Gentle is as follows. He always colors the remaining uncolored vertex in the copy of $K_2$, where the Rascal has just played before.
He chooses his color such that the pair of colors for this $K_2$ is different from all the pairs of colors of the previously totally colored copies of $K_2$. 
Moreover, he ensures that no $K_2$ is monochromatic. Let $\sigma$
 be a nontrivial automorphism of $G$ and let $u\in V(G)$ be  such that $\sigma(u)\neq u$. If $u$ and $\sigma(u)$ belong to the same copy of $K_2$, $\sigma$ does not preserve the coloring. 
Otherwise, the copy of $K_2$ which contains $u$ is sent by $\sigma$ to another one. This implies that  $\sigma$ is not a color preserving automorphism. This shows that the Gentle has a winning strategy with $n+1$ colors. In conclusion, $\GDR(G)=n+1$.
\qed

It is shown in \cite{trenk}, that if $G$ consists of $n$ disjoint copies of $K_2$, then $\displaystyle D(G)=\big\lceil\frac{1+\sqrt{8n+1}}2\big\rceil$. 
Hence, we get examples of graphs for which the $\Rc$-game distinguishing number is quadratic in the classical one. 
The following results are straightforward consequences of Proposition~\ref{prop:worst-inv} and Proposition~\ref{prop:K2}.
\begin{proposition} 
 If $G$ is an involutive graph on $n$ vertices, then:
\begin{itemize}
 \item $\displaystyle D(G)\leq \big\lceil\frac{1+\sqrt{4n+1}} 2\big\rceil$ and
 \item $\displaystyle\GDR(G)\leq \frac{n}{2}+1$. 
\end{itemize}\qed
 \end{proposition}
This result shows that involutive graphs have finite $\Rc$-game distinguishing number. But the bound we obtained does not depend on the classical distinguishing number of the graph. If the Gentle knows this number, having in mind a distinguishing coloring, he can use it to play in a smarter way. 
In other words, for an involutive graph $G$, we can bound the game distinguishing number $\GDR$ using the classical one.
It turns out that $\GDR(G)$ is in this case at most of order $D(G)^2$.
\begin{theorem}\label{theo:inv}
 If $G$ is an involutive graph with $D(G)\geq2$, then $\GDR(G)\leq D(G)^2+D(G)-2$.
\end{theorem}
\proof We set $d=D(G)$. We are going to give a winning strategy for the Gentle with $d^2+d-2$ colors.
Note that, since the Rascal starts the game, the Gentle can play in a way he is never the first to color a vertex in a block.
When the Rascal colors a vertex $u\in V(G)$, the Gentle's strategy is to always answer by coloring the vertex $\overline u$ in the same block. 

\medskip\noindent
Before stating how he chooses the color, we need some definitions.
Let $c$ be a distinguishing coloring of $G$ with $d\geq 2$ colors. We set $k=d^2+d-2$. For all $i,j\in\{1,\cdots,d\}$, we define the following sets of vertices: \mbox{$V_{ij}=\{u\in V(G)|c(u)=i \text{ and } c(\overline u)=j\}$}. 
There are $\frac{d(d-1)}{2}$ sets $V_{ij}$, with $1\leq i<j\leq d$, and $d$ sets $V_{ii}$, with $i\in\{1,...d\}$.
To each set $V_{ij}$, with $1\leq i< j \leq d$, we associate a specific number $\delta_{ij}$ in $\{1,\cdots,\frac{d(d-1)}{2}\}$. 
For each $V_{ij}$, with $i,j\in\{1,...d\}$, we define $r_{ij}$ as follows:
\begin{center}
$r_{ij}=\left\{\begin{array}{ll}
                \delta_{ij}& \text{if } 1\leq i<j\leq d\\
                k-r_{ji} & \text{if } 1\leq j<i\leq d\\
                \frac{d(d-1)}{2}+i &\text{if } i=j \text{ and }1\leq i\leq d-1\\
                0 & \text{if } i=j=d
               \end{array}\right.$
\end{center}
This defines only $d^2$ distinct numbers. However, we need $d^2+d-2$ colors to ensure the second of the following properties.
\begin{enumerate}[{(i)}]
\item  For every $i,i',j,j'\in\{1,\cdots,d\}$, we have $r_{ij}\equiv r_{i'j'} \Mod {k}$ if and only if $i=i' $ and $j=j'$.
\item For every  $i,j\in\{1,\cdots,d\}$, we have $r_{ij}\equiv -r_{ji}\Mod {k}$ if $i\neq j$.
\item For every  $i,j\in\{1,\cdots,d\}$, we have $r_{ii}\not\equiv -r_{jj}\Mod {k}$ if $i\neq j$.
\end{enumerate}
Let now $c': V(G)\rightarrow \{1,\cdots,k\}$ be the coloring built during the game. The Gentle chooses the color as follows. When the Rascal colors $u\in V_{ij}$ with $c'(u)$, he colors $\overline{u}$
such that $c'(\overline u)-c'(u)\equiv r_{ij}\pmod {k}$. 

\medskip\noindent
To prove this strategy is winning, it is enough to show that for any automorphism $\sigma$ which preserves $c'$, we have $\sigma(V_{ij})=V_{ij}$. Since all $V_{ij}$ are monochromatic for the distinguishing coloring $c$, this will show that $\sigma$ also preserves $c$, and $\sigma$
must therefore be the identity.

\medskip\noindent
Let $u$ be a vertex of $G$ and $i,j\in\{1,\cdots,d\}$. Properties (i), (ii), and (iii) easily imply the following facts.
\begin{enumerate}[(a)]
\item If $i\neq j$, then $u\in V_{ij}$ if and only if $c'(\overline u)-c'(u)\equiv r_{ij}\pmod {k}$.
\item $u\in V_{ii}$ if and only if $c'(\overline u)-c'(u)\equiv\pm r_{ii}\pmod {k}$.
\end{enumerate} 
Let $\sigma$ be an automorphism which preserves $c'$. For any vertex $u\in V_{ij}$, with $i,j\in\{1,\cdots,d\}$, we have $c'(\overline u)-c'(u)\equiv c'(\overline{\sigma( u)})-c'(\sigma(u))\pmod {k}$. 
By (a), we get that $c'(\overline{\sigma( u)})-c'(\sigma(u))\equiv r_{ik}\pmod {k}$, and by (b), we have that $\sigma(u)$ is in the same set $V_{ij}$ as $u$.
\qed

In the above proof, decreasing the number of sets $V_{ij}$ tightens the bound. In particular, we have the following corollary.
\begin{corollary}\label{cor:inv}
 Let $G$ be an involutive graph. If there is a coloring with $d\geq 2$ colors, such that ${\rm Bar}$ is the only nontrivial automorphism which preserves this coloring, then
$\GDR(G)\leq 2d-2$.
\end{corollary}
\proof
 We use the same notations as in the proof of Theorem~\ref{theo:inv}. As in this proof, the Gentle uses a coloring $c$ with $d$ colors which is only preserved by ${\rm Bar}$ to build his strategy. 
Since we ask that ${\rm Bar}$ preserves $c$, all the sets $V_{ij}$ with $i\neq j$ are empty. So the Gentle needs only $2d-2$ colors to apply his strategy. Indeed, he plays in a way that for all $u\in V(G)$, 
$c'(\bar u)-c'(u)\equiv i-1\pmod {(2d-2)}$, where $i\in\{1\cdots,d\}$ is such that $u\in V_{ii}$. As for the proof of the previous theorem, we can prove that an automorphism which preserves the coloring $c'$ must also preserve the coloring $c$. Thus, it is either the identity or the ${\rm Bar}$ involution.
All blocks whose vertices do not belong to $V_{11}$ are not monochromatic. 
Hence, the ${\rm Bar}$ involution cannot preserve the coloring $c'$. In conclusion, $c'$ is a distinguishing coloring and the Gentle has a winning strategy with $2d-2$ colors.  
\qed
%
\section{Hypercubes}\label{sec:hypercube}  
In this section, we use Corollary~\ref{cor:inv} to study the game distinguishing number of $Q_n$, the hypercube of dimension $n$ where $n\geq 2$, and prove Theorem \ref{theo:hypercube}. The vertices of $Q_n$ will be denoted by words of length $n$ on the binary alphabet $\{0,1\}$. Let $u(i)$ denote
its $i$-th letter of a vertex $u$. For the classical distinguishing number the question is solved in \cite{bogstad}: $D(Q_3)=D(Q_2)=3$ and $D(Q_n)=2$ for $n\geq 4$. As mentioned before, hypercubes are involutive graphs. The ${\rm Bar}$ involution can be defined as the operation which switches all the letters $0$ to $1$ and all the letters $1$ to $0$. Hence, by Proposition~\ref{prop:order2}, we know that $\GDG(Q_n)=\infty$, for all $n\geq 2$. If we apply directly Theorem~\ref{theo:inv}, we obtain that $\GDR(Q_n)\leq 4$, for $n\geq 2$.
To improve this bound in order to prove the first item of Theorem~\ref{theo:hypercube}, we are going to use Corollary~\ref{cor:inv}. We will have to design nice distinguishing coloring of the hypercube with $2$ colors. 
For this purpose, we introduce determining sets. 

A subset $S$ of vertices of a graph $G$ is  a {\em determining set} if the identity is the only automorphism of $G$ whose restriction to $S$ is the identity on $S$. Determining sets are introduced and studied in \cite{boutin}.
In the following lemma, we give a sufficient condition for a subset of $Q_n$ to be determining.
\begin{lemma}\label{lem:determining}
 Let $S$ be a subset of vertices of $Q_n$, with $n\geq 2$. If for all $i\in\{1,\cdots,n\}$, there are two vertices $v_0^i$ and $v_1^i$ in $S$ such that $v_0^i(i)=0$, $v_1^i(i)=1$ and all the other letters are the same, then $S$ is a determining set of $Q_n$.
\end{lemma}
\proof
Let $\sigma\in {\rm Aut}(G)$, which is the identity on $S$. Suppose for a contradiction that there exists a vertex $u$ in $Q_n\setminus S$ such that $\sigma(u)\neq u$.
Then, there exists $k \in \{1,\cdots,n\}$ such that $u(k)\neq \sigma(u)(k)$. Without loss of generality, we suppose that $u(k)=0$. The vertices $u$  and $v^k_0$ differ on $l$ letters, with $l\in\{1,...n-1\}$.
This means that the distance between $u$ and $v_0^k$ is equal to $l$. Since $\sigma$  fixes $v^k_0$ the distance between $\sigma(u)$  and $v^k_0$ must be $l$, too.
Therefore, the vertices $u$  and $v^k_1$ have $l+1$ distinct letters, whereas the vertices $ \sigma(u)$  and $v^k_1$ have only $l-1$ distinct letters.
This is impossible because the distance between $u$ and $v_1^k$ must be the same as the distance between $\sigma(u)$ and $v_1^k$.
Therefore $\sigma(u)=u$ for all $u\in V(Q_n)$. In conclusion, $S$ is a determining set.\qed

\begin{figure}[h]
\begin{center}
\includegraphics[scale=0.8]{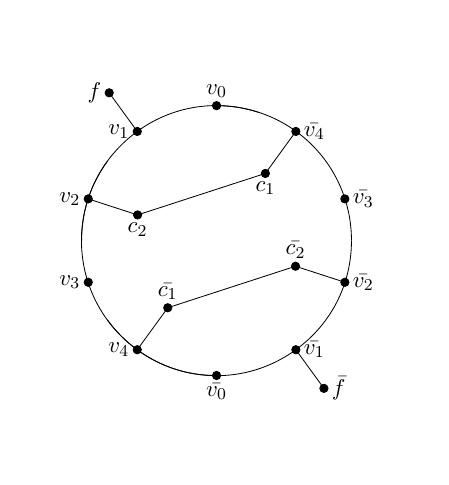}
\caption{\label{fig1} The induced subgraph $S$ in $Q_5$}
\end{center}
\end{figure}

\begin{proposition}
 Let $Q_n$ be the $n$-dimensional hypercube with $n\geq 5$. Then $\GDR(Q_n)=2$.
\end{proposition}
\proof For $i\in\{0,\cdots,n-1\}$, let $v_i$ be the vertex with the $i$ first letters equal to $1$ and the $n-i$ other letters equal to $0$. We define also the following vertices in $Q_n$: $f=10010...0$, $c_1=010...01$ and $c_2=110...01$. The subgraph $S$ will be the graph induced by
$\{v_i,\overline {v_i}|0\leq i<n\}\cup\{f,\overline f, c_1, \overline{c_1}, c_2, \overline{c_2}\}$ (see Figure~\ref{fig1} ). Let $c$ be the coloring with two colors defined by: $c(u)=1$ if 
$u\in V(S)$ and $c(u)=2$ otherwise. We will show that this coloring fits hypothesis of Corollary~\ref{cor:inv}. This would imply that $\GDR(Q_n)=2\times2-2=2$. 

Clearly, ${\rm Bar}$ preserves the coloring $c$ and we prove now this is the only nontrivial color preserving automorphism.
An automorphisms $\sigma$ which preserves this coloring fixes $S$ setwise : $\sigma(S)=S$.  
The restriction of $\sigma$ to $S$, say $\sigma|_S$, is an automorphism of $S$. 
The vertices $v_0, \cdots,v_{n-1}$, and $\overline{v_0}$ show that the hypothesis of Lemma~\ref{lem:determining} holds for the subgraph $S$. 
Hence $\sigma$ is totally determined by the images of the elements of $S$. The two vertices $f$ and $\overline f$ are the only vertices of degree $1$ in $S$.
So, either $\sigma(f)=\overline f$ and $\sigma(\overline f)= f$ or both are fixed. In the first case, this implies that $\sigma(v_1)=\overline{v_1}$ 
and $\sigma(\overline{v_1})=v_1$. In the second case, $v_1$ and $\overline{v_1}$ are also fixed by $\sigma$. This is
because they are respectively the only neighbors of $f$ and $\overline f$ in $S$. The vertices $v_0$ and $v_2$ are the two remaining neighbors of $v_1$.
Since they have not the same degree in $S$, they cannot be switched by $\sigma$. 
Hence, in the first case, $\sigma(v_i)=\overline{v_i}$ and $\sigma(\overline{v_i})=v_i$, with $i\in\{0,2\}$. In the second case, $v_i$ and $\overline{v_i}$, with $i\in\{0,2\}$, are fixed by $\sigma$. 
There is exactly one path of size $n-2$ between $v_0$ and $\overline{v_2}$ (resp.  $v_2$ and $\overline{v_0}$). Hence, in the first case these paths are switched, whereas they are fixed pointwise in the second.  
After that, it is easy to show that $\sigma|_S$ is either the identity or the ${\rm Bar}$ involution.
\qed

\begin{figure}[h]
\begin{center}
\includegraphics[scale=0.8]{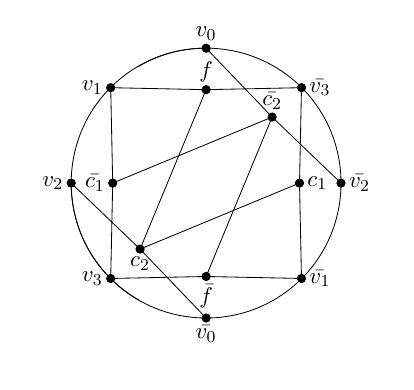}
\caption{The induced subgraph $S$ in $Q_4$}\label{fig2}
\end{center}
\end{figure}

The above proof fails for $Q_4$ because the subgraph $S$ will have automorphisms other than the {\rm Bar} involution (see Figure \ref{fig2}). Nevertheless, we can decrease the upper bound from 4 to 3. 
\begin{proposition}
 We have $\GDR(Q_4)\leq 3$.
\end{proposition}
\proof Let $S$ be the subgraph of $Q_4$ induced by the five vertices $0000$, $1000$, $1100$, $1110$ and $1011$. Let $\overline S$ be the subgraph induced by the opposite vertices
of those in $S$. Note that $S$ and $\overline S$ are disjoint. The Gentle's strategy is as follows. When the Rascal colors a vertex $u$ in $V(G)$, the Gentle
colors $\overline u$. He chooses the coloring $c$ according to these rules:
\begin{itemize}
 \item if $u$ is in $V(S)$, then $c(\overline u)-c(u)\equiv 1 \pmod {3}$, 
 \item if $u$ is in $V(\overline S)$, then $c(\overline u)-c(u)\equiv -1\pmod {3}$,
 \item $c(\overline u)-c(u)\equiv 0\pmod {3}$ otherwise.
\end{itemize}
Let $\sigma$ be a color preserving automorphism. Since $c(\bar u)-c(u)\equiv 1\pmod{3}$ if and only if $u\in S$, we have that $S$ is fixed setwise by $\sigma$. Let $\sigma|_S$ be the restriction of $\sigma$ to $S$. The vertices $0000$ and $1110$ are the only
vertices of degree $1$ in $S$. Then either $\sigma|_S$ fixes them or switches them. But in $Q_4$, we have $d(0000,1011)\neq d(1110,1011)$. Since $\sigma|_S$ comes from an automorphism of $Q_n$, it must preserve distances in $Q_n$. 
But, the vertex $1011$ is fixed by $\sigma|_S$ because it is the only isolated vertex of $S$. Hence, we have that $\sigma|_S$ fixes $0000$ and $1110$. Moreover, it is clear that all the vertices of $S$ are fixed by $\sigma|_S$ and hence by $\sigma$. The vertex
$1111$ is also fixed  by $\sigma$, because it is the opposite vertex of $0000$. Finally, all the vertices of $V(S)\cup\{1111\}$ are fixed by $\sigma$ and this subset fits hypothesis of Lemma~\ref{lem:determining}. This shows that $\sigma$ is the identity of $Q_4$.   
\qed

To complete the proof of Theorem~\ref{theo:hypercube}, it remains to settle the case of hypercubes of dimension $2$ and $3$.
For the square $Q_2$, isomorphic to $C_4$, an easy computation shows that $\GDR(Q_2)=3$. For $Q_3$, it turns out that $\GDR(Q_3)$ is also equal to $3$.
It comes from the fact that the complementary graph of $Q_3$ is isomorphic to $K_4\square K_2$ and from the following proposition.

We will write the vertices of $K_4\square K_2$ as couples $(i,x)$, with $i\in \{1,\cdots,4\}$ and $x\in\{l,r\}$. The first coordinate is the one
associated to $K_4$ and the second the one associated to $K_2$. We denote by $K_2^i$ the $K_2$-fiber whose two vertices have $i\in\{1,\cdots,4\}$ as first coordinate.
Let $c$ be a coloring of $K_4\square K_2$. We say that two distinct $K_2$-fibers, $K_2^i$ and $K_2^j$, with $i,j\in\{1,\cdots,4\}$ are {\em colored the same} if $c(K_2^i)=c(K_2^j)$.
Moreover, if $c((i,l))=c((j,l))$, we say that the two fibers are {\em strictly colored the same}. 
\begin{proposition}
 We have $\GDR(K_4\square K_2)=3$.
\end{proposition}
\proof
Since $D(K_4\square K_2)=3$, we only have to prove that the Gentle has a winning strategy with $3$ colors.
The Rascal starts by coloring the vertex $(1,l)$ with $1$. The Gentle replies by coloring with $2$ the vertex
$(2,l)$. Now, there are two cases.

\noindent
\textbf{Case 1:} The Rascal colors the second vertex of $K_2^1$ or $K_2^2$ with a colors $\delta\in\{1,2,3\}$.
Without loss of generality we assume he colors $(1,r)$. There are two subcases.

\textbf{Subcase 1.1:} $\delta\in\{1,3\}$. The Gentle colors $(2,r)$ with $2$. The Gentle's strategy is now to play in the same $K_2$-fiber as Rascal. Since the two $K_2$-fibers $K_2^1$ and $K_2^2$ share no color, the Gentle can ensure that at the end of the game no pair of $K_2$-fibers 
are colored the same and $K_2^3$ is not monochromatic. This shows that a color preserving automorphism $\sigma$ cannot switch the $K_2$-fibers. Hence 
$\sigma(K_2^3)=K_2^3$ and, since this fiber is not monochromatic, we have that $\sigma$ cannot switch the $K_4$-fibers, too. Finally, $\sigma$ is the identity.

\textbf{Subcase 1.2:} $\delta=2$. The Gentle colors $(2,r)$ with color $1$.
By coloring the vertices in the same $K_2$-fiber as the vertices the Rascal will play, the Gentle can ensure that at the end of the game
$K_2^1$ and $K_2^2$ are the only $K_2$-fibers colored the same. Moreover, he can ensure that one $K_2$-fiber, say $K_2^3$,  is not monochromatic. 
Let $\sigma$ be a color preserving automorphism. 
The fibers $K_2^3$ and $K_2^4$ are fixed setwise by $\sigma$. Since $K_2^3$ is not monochromatic, $\sigma$ cannot switch the $K_4$-fibers. 
We conclude that every $K_2$-fibers are fixed pointwise by $\sigma$, which means $\sigma$ is trivial. 

\noindent
\textbf{Case 2:} The Rascal colors the vertex $(3,x)$, with $x\in\{l,r\}$. The Gentle answers by coloring the vertex $(4,x)$.
The Gentle chooses his colors such that we can assume $(1,l)$ and $(3,x)$ are colored with $1$, $(2,l)$ is colored with $2$ and $(4,x)$ is colored with $3$ (maybe we need to permute colors $1$ and $2$).
We only deal with the case where $x=l$. The case where $x=r$ is very similar.
When the Rascal's move is to color $(1,r)$ (resp. $(3,r)$), the Gentle answers by coloring $(3,r)$ (resp. $(1,r)$).
When the Rascal's move is to color $(2,r)$ (resp. $(4,r)$), the Gentle answers by coloring $(4,r)$ (resp. $(2,r)$).
He can choose the colors he uses, such that no pair of $K_2$-fibers is strictly colored the same. Moreover, one of the four $K_2$-fibers is not monochromatic and
is not colored the same as each of the three other $K_2$-fibers. At the end of the game, we are in the same situation as in subcase \textbf{1.2}. In conclusion, 
the coloring built in this case is also 3-distinguishing.
\qed

\section{Cycles}\label{sec:cycle}
\begin{figure}[h]
\begin{center}
\includegraphics[scale=0.9]{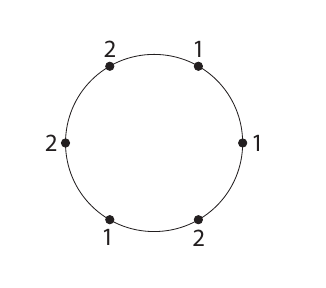}
\caption{The unique distinguishing $2$-coloring of $C_6$}\label{fig0}
\end{center}
\end{figure}

 In this section, we prove Theorem \ref{theo:cycle} about cycles. We start with even cycles which are involutive graphs. For such graphs, the opposite of a vertex can be defined as the unique vertex at maximal distance.
Corollary~\ref{cor:inv} enables us to conclude when the size of the even cycle is greater or equal to $12$. 
For $C_4$, it is easy to verify that $\GDR(C_4)=3$. Using the fact that the coloring shown in Figure~\ref{fig0} is the only distinguishing $2$-coloring of $C_6$, it is easy to see that $\GDR(C_6)=3$. For $C_8$ and $C_{10}$, the Gentle's strategies are less obvious. We give them in the below proposition.
\begin{proposition} 
We have $\GDR(C_8)=\GDR(C_{10})=2$. 
\end{proposition}
\proof
We start with Gentle's winning strategy for $C_8$.
Let $(x_1,\cdots,x_8)$ be the cycle of order $8$. Let $c$ be a $2$-coloring with exactly three vertices colored with color $2$. If these three vertices do not induce
a stable set or a $P_3$, then $c$ is $2$-distinguishing. The Rascal starts by coloring $x_1$ with color $1$. The Gentle replies by coloring
$\overline{x_1}$ with the same color. Gentle's strategy is now to use the other color than the one the Rascal played just before.
In this way, the number of vertices colored with  color $2$ will be exactly three. To avoid that
these three vertices induce a stable set or a $P_3$, he plays exactly one move in each of the pairs $\{x_2,x_4\},\{x_6,x_8\},\{x_3,x_7\}$.

Now, we give Gentle's strategy for $C_{10}$. Let $(x_1,\cdots,x_{10})$ be the cycle of order $10$. Rascal's first move is to color $x_1$ with color $1$. The Gentle replies by coloring
$\overline{x_1}$ with the same color. After that there are two cases.

\noindent
\textbf{Case 1:} The Rascal uses color $1$ for his second move. Using symmetries, we obtain only two subcases. In both subcases, the Gentle answers first with color $1$. 
Then, he uses the other color than the one used by the Rascal just before. 

\textbf{Subcase 1.1:} The Rascal colors $x_2$ with $1$. The Gentle colors $x_7$ with $1$. Then, he plays using the pairs $\{x_3,x_5\},\{x_8,x_{10}\},\{x_4,x_7\}$. 

\textbf{Subcase 1.2:} The Rascal colors $x_3$ with $1$. The Gentle answers by coloring $x_4$ with $1$. Then, he plays following the pairs $\{x_2,x_5\},\{x_7,x_9\},\{x_8,x_{10}\}$.
%

\medskip\noindent
\textbf{Case 2:} The Rascal use color $2$ for his second move. The Gentle answers by coloring a vertex with color $2$, too. By symmetry, we assume that  $x_2$
and $x_3$ are colored with $2$. The Gentle plays according to the pairs $s_1=\{x_4,x_5\}$, $d=\{x_7,x_8\}$ and $s_2=\{x_9,x_{10}\}$. If the Rascal plays in $s_1$ or $s_2$,
 the Gentle copies the color he used, whereas, when the Rascal plays in $d$, the Gentle plays the other color. 
 
 We let the reader check that the strategies for $C_8$ and $C_{10}$ are indeed winning. 
\qed 
\begin{proposition} 
For $n\geq 6$, we have $\GDR(C_{2n})=2$.
\end{proposition}
\proof Let $u,v,w$ be vertices of $C_{2n}$ such that $d(u,v)=1$, $d(u,w)=3$ and $d(v,w)=2$. We set $S=\{u,\bar u, v, \bar v, w ,\bar w\}$ and define the coloring $c$ with $2$ colors as follows: $c(x)=1$ if $x\in S$ and $c(x)=0$ otherwise. 
Let $\sigma$ be a nontrivial color preserving automorphism. The subgroup $<\sigma>$ acts on the set of the three blocks contained in $S$. If $\sigma$ is an axial reflection then one block must be stable under $\sigma$. This block is the axis of the reflection. Since $2n\geq 12$, we easily check that the two other blocks cannot be map to each other by this reflection. Since $\sigma$ cannot be a reflection, $\sigma$ must be a rotation. The only possibility for a rotation to let $S$ stable is to map each vertex to its opposite. 
Hence $\sigma$ restricted to $S$ is the ${\rm Bar}$ involution. But it is straightforward to show that $\sigma$ must be the ${\rm Bar}$ involution on the whole cycle. This shows that the coloring $c$ fits Corollary~\ref{cor:inv} hypothesis. We conclude that $\GDR(C_{2n})\leq 2$.\qed

We compute now the $\Gc$-game distinguishing number for odd, but not prime cycles. Odd cycles are not involutive graphs, but when the number $n$ of vertices is not prime, the action of the automorphism group is not primitive. 
We are going to use a complete block system, with blocks of size equal to the least prime divisor of $n$.   



\begin{proposition}\label{prop:odd-cycle}
 If $C_n$ is an odd cycle such that $n>9$  and $n$ is not prime, then $\GDG(C_n)=2$.
\end{proposition} 
\proof
Let $p$ be the least prime divisor of $n$. Since $n\geq 15$, we have $n=kp$, with $k\geq5$.
Let $(x_1,\cdots,x_n)$ be the cycle $C_n$. We define $k$ disjoint subsets of vertices $V_j=\{x_{j+lk}|l\in\{0,\cdots,p-1\}\}$, with $j\in\{1,\cdots,k\}$.
The set $V_1$ will play a particular role. We denote by $\Delta_i$, with $i\in\{1,\cdots,p\}$, the reflection with fixed point $x_{1+(i-1)k}$.
Note that the automorphisms $\Delta_i$, with $i\in\{1,\cdots,p\}$, moves the sets $V_j$, with $j\in\{1,\cdots,k\}$, in the following way:
$\Delta_i(V_1)=V_1$ and for $j\neq 1$, $\Delta_i(V_{j})=V_{k+2-j}$. Since $k$ is odd, $j\neq k+2-j$ and $\Delta_i(V_j)\neq V_j$, if $j\neq 1$.

We describe the Gentle's strategy. He starts by coloring $x_1$ with color $1$. 
More precisely, he chooses the fixed point of $\Delta_1$.
\begin{itemize}
 \item If the Rascal colors a vertex $x$ in $V_1$, he answers by coloring the vertex $\Delta_1(x)$ in $V_1$ with the other color. Note that since the Gentle starts the game
by coloring the fixed point of $\Delta_1$, this is always possible. Therefore, $\Delta_1$ is broken and there will be exactly $1+\frac{p-1}2$ vertices colored
with $1$ in $V_1$ at the end of the game.
 \item If the Rascal plays in $V_j$, with $j\neq 1$, then he plays in $V_{k+2-j}$. Except if the Rascal just colors the penultimate vertex of $V_j$. In that case he plays
in $V_j$ too, and colors the last vertex in a way the parity of the number of $1$ in $V_j$ is not the same as the parity of $1+\frac{p-1}{2}$.
The vertex he chooses to color in $V_{k+2-j}$ is determined as follows.
If at least one reflection $\Delta_i$, with $i\in\{2,\cdots,p\}$, is not already broken by the coloring, he chooses one of them, say $\Delta_{i_0}$.
We call $V_j'$ the set of already colored vertices in $V_j$. There is exactly one more colored vertex in $V_j$ than in $V_{k+2-j}$. Hence $\Delta_{i_0}(V_j')$
 has at least an uncolored vertex, say $u$. The Gentle colors $u$ such as $u$ and $\Delta_{i_0}^{-1}(u)$ are not of the same color.
If all the symmetries $\Delta_i$ are already
broken, he plays randomly in $V_{k+2-j}$.
\end{itemize}
At the end of the game, there are exactly $1+\frac{p-1}{2}$ vertices colored with $1$ in $V_1$.
It is clear that Gentle's strategy leads him to play the last vertex of any set $V_j$. Hence, in all sets $V_j$, with $j\neq 1$, the number of $1$s has not the same
parity as in $V_1$.
This proves that the only possible color preserving automorphisms are those under which $V_1$ is stable. The Gentle plays $\frac{(k-1)p}{2}$ moves outside $V_1$.
There are only $k-1$ moves used to control the parity of $1$s. The number of remaining moves to break the $p-1$ symmetries $\Delta_2, \cdots,\Delta_{p}$ is
$\frac{(p-2)(k-1)}{2}$. Since $k\geq p$ and $k\geq 5$, this number is greater or equal to $p-1$.
This shows that all the symmetries under which $V_1$ is stable are broken.
A color preserving automorphism $\sigma$ is finally either a rotation under which $V_1$ is stable or the identity.
Assume $\sigma$ is a rotation. The subgroup $<\sigma>$ acts naturally on $V_1$. Since $|V_1|$ is prime, $<\sigma>$ must be of order $p$ and acts transitively on $V_1$.
This is impossible because $V_1$ is not monochromatic. In conclusion, $\sigma$ must be the identity and the Gentle has a winning strategy with two colors.
\qed

For $C_9$, the general strategy fails but two colors are still enough for Gentle to win.

\begin{proposition}
 We have $\GDG(C_9)=2$.
\end{proposition}
\proof
Let $(x_1,\cdots,x_9)$ be the cycle of order $9$ and $c$ the coloring built throughout the game. 
The Gentle can play such that after his second move, the game is in one of the following cases.

\noindent\textbf{Case 1:} $c(x_1)=c(x_4)=1$ and $c(x_7)=2$.
The Gentle applies the same strategy as in Proposition~\ref{prop:odd-cycle}. 

For the other cases, the Gentle defines three pairs $S,O,O'$ of uncolored vertices. He plays in a way the pairs $O$ and $O'$ have vertices of different colors, and the pair $S$ is monochromatic.

\noindent\textbf{Case 2:}  $c(x_1)=c(x_2)=1$ and $c(x_6)=2$. The pairs are defined by $O=\{x_3,x_9\}$, $O'=\{x_5,x_7\}$  and $S=\{x_4,x_8\}$. 

\noindent\textbf{Case 3:} $c(x_1)=c(x_3)=1$ and $c(x_2)=2$. He defines $S=\{x_4,x_6\}$, $O=\{x_5,x_8\}$  and $O'=\{x_7,x_9\}$. 

\noindent\textbf{Case 4:}  $c(x_1)=c(x_5)=1$ and $c(x_3)=2$. He makes $S=\{x_2,x_4\}$, $O=\{x_6,x_9\}$  and $O'=\{x_7,x_8\}$. 

\qed

To finish the proof of Theorem~\ref{theo:cycle}, we deal with prime cycles. Note that, in that case, the automorphism group action is primitive. Once more, the small case $C_5$ does not fit the general strategy. By computation, we get $\GDG(C_5)=3$.  
\begin{proposition}\label{propo:primecycle}
 If $C_p$ is a cycle of prime order $p>5$, then $\GDG(C_p)\leq 3$.
\end{proposition}
 \proof Assume that three colors are allowed during the game. We begin the proof by showing the Rascal cannot win if he does not always play the same color as the one used by the Gentle just before.
We denote respectively by $n_1$, $n_2$ and $n_3$,
the number of vertices colored with $1$, $2$ and $3$ during the game.
 As long as the Rascal  copies the color played by the Gentle just before, the three numbers $n_1$, $n_2$ and $n_3$ are even at the end of the Rascal's turn. The first time the Rascal plays a different color,
two of the three numbers, $n_1$, $n_2$ and $n_3$ are odd at the end of his turn. Without loss of generality, we can say $n_1$ and $n_2$ are odd. Then the Gentle colors a vertex with
$3$. The numbers $n_1$, $n_2$ and $n_3$ are now all odd. Until the end of the game, the Gentle strategy is now to play the same color as the one the Rascal just used before him. In this way,
 he preserves the parity of $n_1$, $n_2$ and $n_3$. At the end of the game, the coloring is such that $n_1$, $n_2$ and $n_3$ are odd.
We show that this coloring is a distinguishing coloring.
Let $\sigma$ be a color preserving automorphism. The size of an orbit $\Oc$ under the action of $<\sigma>$ must divide $|<\sigma>|$.
But the cardinality of $<\sigma>$ is either $1$, $2$ or $p$.
Then, this orbit has $1$, $2$ or $p$ as its size.
None of $n_1$, $n_2$ and $n_3$ are null, so the coloring is not monochromatic and then the size of $\Oc$ cannot be $p$. The automorphism $\sigma$ cannot be a rotation
because, for prime cycle, a rotation acts transitively and has an orbit of size $p$. Assume $\sigma$ is a reflection. Then, there is exactly one fixed point. 
Without loss of generality, we suppose it is colored $1$.
All the vertices colored with $2$ are in orbits of size $2$. Since $n_2$ is not null, there are $k>0$ such orbits.
We get $n_2=2k$, which is a contradiction. Hence $\sigma$ must be the identity.

We can now assume that the Rascal always copies the color played by the Gentle just before. Playing a different move will actually lead him to defeat, whatever has happened before. 
The winning strategy for the Gentle is as follows.
He starts with $1$ and his second move is to color a vertex with $2$. The Rascal's second move must also be $2$. Let $\sigma_2$ be the reflection which switches the two
vertices colored with $2$ and $v$ the vertex fixed by $\sigma_2$. If $v$ is already colored, the Gentle plays $1$ wherever he wants. If not, he colors $v$ with $1$.
After that, the Gentle always uses the color $1$, except for his last turn, where he will play $3$. Since $p>5$, he really has enough turns to play this $3$.
At the end of the game, the coloring has the following properties: exactly two vertices are colored with $2$ and exactly one vertex is colored with $3$.
Since there is only one vertex colored with $3$ no rotation can preserve the coloring. Since there are only two vertices colored with $2$, the only possible colors
preserving reflection is $\sigma_2$. But the Gentle played in a manner that the unique $3$ is not the fixed point of this reflection. In conclusion, the coloring is distinguishing and the Gentle wins with $3$ colors.
\qed 
\section{Conclusion and further works}
In this article, we have defined two new game invariants of graphs, in the same spirit as the game chromatic number or as the game domination number. 
Since these invariants could be infinite, we started by giving some sufficient conditions for a graph to have infinite game distinguishing numbers.
But a total characterization seems far to be found and sounds like a very challenging open problem. We propose the following conjecture.
\begin{conjecture}
 If $G$ is a graph with no automorphism of order $2$, then $\GDG(G)$ and $\GDR(G)$ are both finite.
\end{conjecture}

After that, we defined a new class of graphs, the involutive graphs, for which we can bound quadratically the $\Gc$-game distinguishing number using informations on the classical one. 
Results on this class are then applied to compute the exact value of the game distinguishing number of hypercubes and even cycles. For the hypercube of dimension $4$, determining if $\GDR(Q_4)$ is equal to $2$ or $3$ remains open. 

For odd but not prime cycles, we were also able to compute the exact value of the two game invariants. The Gentle strategy made an intensive use of the imprimitivity of the automorphism group action. 
When the cycle is prime, we only know that $\GDG$ is bounded above by $3$. This bound is sharp for small cases. Indeed, $\GDG(C_5)=\GDG(C_7)=3$. But, it seems this is more related to the small size of these cycles than to the primality of their order. Supported by computer experimentation, we state the following conjecture.
\begin{conjecture}
 If $C_p$ is a cycle of prime order $p\geq 11$, then $\GDG(C_p)=2$.
\end{conjecture}

Remark that for prime cycles, the automorphism group action is primitive. It would be interesting to study other graphs for which it is the case. It seems that the question is harder than in the imprimitive case. 

\section*{Acknowledgements}

The research was financed by the ANR-14-CE25-0006 project of the French National Research Agency.
\nocite{*}
\bibliographystyle{}

\end{document}